\documentclass[11pt]{article}
\usepackage[T1]{fontenc}
\usepackage{lmodern}
\usepackage[hmargin=1in,vmargin=0.93in]{geometry}
\usepackage{microtype}
\usepackage{amsmath,amssymb,amsthm}
\usepackage{xcolor}
\usepackage[colorlinks=true,allcolors=blue!55!black]{hyperref}
\hypersetup{pdftitle={A disproof of the El-Zahar--Erd\H{o}s conjecture},pdfauthor={Chenghua Liu, Boning Meng}}

\newtheorem{theorem}{Theorem}
\newtheorem{lemma}[theorem]{Lemma}
\newcommand{\ee}{\mathrm e}

\title{A disproof of the El-Zahar--Erd\H{o}s conjecture}
\author{%
  Chenghua Liu, Boning Meng\\[0.25em]
  \small Renormalization Group\enspace\textperiodcentered\enspace
    \href{https://r13n.group/}{r13n.group}\\
  \small Correspondence:
    \href{mailto:renormalization@r13n.group}{renormalization@r13n.group}%
}
\date{}

\begin{document}
\maketitle
\vspace{-2.6em}

\begin{abstract}
El-Zahar and Erd\H{o}s asked whether every graph of bounded clique
number and sufficiently large chromatic number contains two
vertex-disjoint induced subgraphs, each of large chromatic number,
with no edge between them. We give a negative answer even for
triangle-free graphs. For every integer \(g\ge 4\), we show by a
probabilistic argument that there exist graphs of girth at least \(g\)
and arbitrarily large chromatic number in which any two
vertex-disjoint induced subgraphs that each require at least seven
colors are joined by an edge.
\end{abstract}

\section{Introduction}

All graphs are finite and simple. Write \(\omega(G)\) for the clique
number of \(G\) and \(\chi_G(X)=\chi(G[X])\) for the chromatic
number of the subgraph induced by \(X\subseteq V(G)\), with
\(\chi_G(\varnothing)=0\). Two vertex sets are \emph{anticomplete}
if they are disjoint and no edge joins them.
The \emph{girth} of a graph is the length of its shortest cycle,
with girth \(+\infty\) for a forest.

El-Zahar and Erd\H{o}s \cite{ElZaharErdos} proved that, for all positive
integers \(t,c\), every graph \(G\) with \(\omega(G)<t\)
and sufficiently large chromatic number contains anticomplete sets
\(A,B\) with \(\chi_G(A)\ge 3\) and \(\chi_G(B)\ge c\). They asked
whether one can also require \(\chi_G(A)\ge c\), allowing the threshold
for \(\chi(G)\) to depend on \(t\) and \(c\). Nguyen, Scott, and Seymour
\cite{NSS} showed that one can require \(G[A]\) to have minimum degree
at least \(c\) while keeping \(\chi_G(B)\ge c\).

\begin{theorem}
\label{thm:main}
For every pair of integers \(d\ge 1\) and \(g\ge 4\), there exists a
graph \(G\) with \(\chi(G)\ge d\) and girth at least \(g\) such that
every anticomplete pair \(A,B\) in \(G\) satisfies
\[
 \min\bigl(\chi_G(A),\chi_G(B)\bigr)\le 6.
\]
\end{theorem}

Taking \(g=4\) gives triangle-free graphs, so \(\omega(G)<3\).
These graphs have arbitrarily large chromatic number, yet no
anticomplete pair has both sides of chromatic number at least \(7\).
Thus the answer to the question of El-Zahar and Erd\H{o}s is negative
already for \(t=3\) and \(c=7\).

We construct a random graph \(H\) on blocks \(V_1,\ldots,V_r\) of
equal size, adding edges only between distinct blocks. The edge
probabilities increase sharply along a prescribed order of the block
pairs. Each pair \((V_i,V_j)\) has a size threshold, and we say that
\(B\) is \emph{large on} the pair if both \(|B\cap V_i|\) and
\(|B\cap V_j|\) exceed the threshold for that pair.

The ordering ensures that whenever \(A\) and \(B\) are anticomplete,
one of them, say \(B\), has the following property. We can cut the
sequence \(V_1,\ldots,V_r\) after some \(V_s\) so that every pair on
which \(B\) is large has one block on each side of the cut. Thus every
edge of \(H[B]\) lying entirely on one side comes from a pair on which
\(B\) is not large. A sparsity estimate shows that the vertices of
\(B\) on each side induce a \(2\)-degenerate graph and hence can be
colored with three colors~\cite[Proposition~5.2.2]{Diestel}. Using
disjoint sets of three colors on the two sides gives a \(6\)-coloring
of \(H[B]\). Finally, we delete one vertex from each cycle of length
less than \(g\). The resulting induced subgraph \(G\) has girth at
least \(g\). We show that there are few enough such cycles for \(G\)
to still have chromatic number at least \(d\).

\section{An ordered block construction}

Let \(V_1,\ldots,V_r\) be independent sets of size \(n\) partitioning
\(V(H)\). Assign a parameter \(D_{ij}>1\) to each pair \((i,j)\)
with \(1\le i<j\le r\). For \(X\subseteq V(H)\), write \(X_i=X\cap V_i\).
Call \((i,j)\) \emph{weak for \(X\)} if
\(\min\{|X_i|,|X_j|\}\le nD_{ij}^{-8}\). Let \(W_X\) be the spanning
subgraph of \(H[X]\) consisting of all edges joining \(X_i\) to \(X_j\)
for pairs \((i,j)\) weak for \(X\).

Before sampling any edges, we fix the following order \(\prec\) on the
block pairs. We order pairs by increasing second index, breaking ties
by decreasing first index. For three blocks \(V_1,V_2,V_3\), this gives
\[
 (1,2)\prec(2,3)\prec(1,3).
\]
We choose the parameters \(D_{ij}\) to satisfy \textup{(O)} and require
the graph \(H\) to satisfy \textup{(L)} and \textup{(S)}.
\begin{enumerate}
\item[(O)] For any two block pairs \((i,j)\prec(k,\ell)\), we have
\(D_{k\ell}>D_{ij}^{24}\).
\item[(L)] For each \(i<j\), if \(U\subseteq V_i\) and \(V\subseteq V_j\)
both have size at least \(nD_{ij}^{-1/3}\), then an edge joins \(U\) to \(V\).
\item[(S)] For every nonempty \(X\subseteq V(H)\), the graph \(W_X\)
has fewer than \(3|X|/2\) edges.
\end{enumerate}

By \textup{(O)}, the size required by \textup{(L)} on a later pair is
smaller than the cutoff for weakness on any earlier pair. Indeed,
\begin{equation}
 nD_{k\ell}^{-1/3}<nD_{ij}^{-8}
 \qquad\text{whenever } (i,j)\prec(k,\ell).
 \label{eq:scales}
\end{equation}

\begin{lemma}
\label{lem:order}
Suppose that \textup{(O)}, \textup{(L)}, and \textup{(S)} hold.
Then every anticomplete pair \(A,B\) in \(H\) has a side of
chromatic number at most \(6\).
\end{lemma}

\begin{proof}
Fix \(X\subseteq V(H)\). A pair weak for \(X\) remains weak for every
\(Y\subseteq X\), so \(W_X[Y]\) is a subgraph of \(W_Y\). By \textup{(S)},
every nonempty induced subgraph of \(W_X\) has average degree less
than three and hence a vertex of degree at most two. Thus \(W_X\) is
\(2\)-degenerate and admits a \(3\)-coloring.

Let \(F_X\) be the graph on the block indices \(1,\ldots,r\) whose
edges are the pairs that are not weak for \(X\). A proper coloring of
\(F_X\) assigns a color to each block and separates the endpoints of
every edge of \(H[X]\) outside \(W_X\). Pairing this block color with
the vertex color in \(W_X\) gives
\begin{equation}
 \chi_H(X)\le 3\chi(F_X).
 \label{eq:product-color}
\end{equation}
By \eqref{eq:product-color}, it suffices to show that \(F_A\) or \(F_B\) is bipartite.
This holds if either graph is edgeless. Otherwise, exchange \(A\) and
\(B\) if necessary so that the earliest edge \((i,j)\) in
\(E(F_A)\cup E(F_B)\) belongs to \(F_A\).

We claim that every edge \((k,\ell)\) of \(F_B\) satisfies
\(k\le i<j\le\ell\). By the choice of \((i,j)\), we have \(\ell\ge j\),
and if \(\ell=j\), the ordering within that second index gives \(k\le i\).
Suppose then that \(\ell>j\) and \(k>i\). The pair \((i,\ell)\) follows
both \((i,j)\) and \((k,\ell)\), so \eqref{eq:scales} gives
\[
 nD_{i\ell}^{-1/3}<\min\{nD_{ij}^{-8},nD_{k\ell}^{-8}\}.
\]
Since \((i,j)\in E(F_A)\) and \((k,\ell)\in E(F_B)\), the two
quantities on the right are smaller than \(|A_i|\) and \(|B_\ell|\),
respectively. Property \textup{(L)} therefore gives an edge from
\(A_i\) to \(B_\ell\), contradicting anticompleteness.
This proves the claim. Every edge of \(F_B\) crosses the cut after
index \(i\), so \(F_B\) is bipartite.
\end{proof}

\section{Random blocks and removal of short cycles}

To prove existence, we show that with probability tending to one as the
block size grows, a single random graph \(H\) satisfies both \textup{(L)}
and \textup{(S)} and has few short cycles. For \textup{(L)}, two large
sets in distinct blocks give many possible edges between them, and it
suffices for just one to be present. For \textup{(S)}, consider a block
pair that is weak for \(X\). One block contains so few vertices of
\(X\) that each vertex of \(X\) in the other block has few neighbors
among them in expectation. The proof turns this observation into a bound
on the total number of edges from weak pairs, valid simultaneously
for all \(X\).

\begin{lemma}
\label{lem:random}
For each integer \(r\ge2\) there is a constant \(D_0(r)\) with the
following property. Fix an integer \(g\ge4\) and parameters
\(D_{ij}\ge D_0(r)\) for \(1\le i<j\le r\). On \(r\) independent
blocks of size \(n>\max_{i<j} D_{ij}\), form a random graph \(H\)
by including each possible edge between \(V_i\) and \(V_j\)
independently with probability \(D_{ij}/n\). With probability tending
to one as \(n\to\infty\), \(H\) satisfies \textup{(L)} and \textup{(S)}
and contains at most \(\sqrt n\) cycles of length less than \(g\).
\end{lemma}

\begin{proof}
Write \(D_*=\min_{i<j}D_{ij}\) and \(N=rn\). Fix \(i<j\) and put
\(D=D_{ij}\). For subsets \(U\subseteq V_i\) and \(V\subseteq V_j\)
of size at least \(nD^{-1/3}\), the probability that no edge joins
them is at most \(\exp(-D|U||V|/n)\le\exp(-nD^{1/3})\).
There are at most \(4^n\) choices of \(U,V\) per block pair, so the
probability that \textup{(L)} fails is at most
\(\binom r2\exp\bigl(n(\log 4-D_*^{1/3})\bigr)=o(1)\) whenever
\(D_*^{1/3}>\log 4\).

For \textup{(S)}, fix a set \(X\) of size \(m\ge1\) and write
\(Z_X=|E(W_X)|\). We retain the dependence on \(m\) in our estimate
so that we can sum the resulting probability bound over all
\(\binom Nm\) choices of \(X\). Consider a pair weak for \(X\), with parameter \(D\)
and intersection sizes \(u\le v\). Its contribution to \(\mathbb E Z_X\)
is \(vDu/n\). Using \(u/n\le D^{-8}\) and \(u\le m\), we obtain
\[
 \frac{Du}{n}
 =D\left(\frac{u}{n}\right)^{1/4}
    \left(\frac{u}{n}\right)^{3/4}
 \le D^{-1}\left(\frac{m}{n}\right)^{3/4}.
\]
The sum of the larger sizes \(v\) over all weak pairs is at most
\((r-1)m\), since each block occurs in at most \(r-1\) pairs. Hence
\[
 \mu_X:=\mathbb E Z_X
 \le \frac{rm}{D_*}\left(\frac{m}{n}\right)^{3/4}
 = \frac{r^{7/4}m}{D_*}\left(\frac{m}{N}\right)^{3/4}.
\]
Choose \(D_0(r)\ge 4\ee r^{7/4}\), so the preceding bound also gives \(\ee\mu_X/m\le (m/N)^{3/4}/4\).
For fixed \(X\), its weak pairs are determined before any edges are
sampled. Thus \(Z_X\) is a sum of independent indicators, and for
\(q=\lceil3m/2\rceil\) we have
\[
 \Pr(Z_X\ge q)\le\mathbb E\binom{Z_X}{q}
 \le\frac{\mu_X^q}{q!}
 \le\left(\frac{\ee\mu_X}{q}\right)^q
 \le\left[\frac14\left(\frac{m}{N}\right)^{3/4}\right]^{3m/2}.
\]
Multiplying by \(\binom Nm\le(\ee N/m)^m\) leaves a factor of
\((m/N)^{9m/8-m}=(m/N)^{m/8}\). Thus the probability that some
\(m\)-vertex set violates \textup{(S)} is at most
\[
 \left[\frac{\ee}{8}\left(\frac{m}{N}\right)^{1/8}\right]^m.
\]
For \(m\le\sqrt N\), these bounds sum to at most
\(\sum_{m\ge1}(N^{-1/16}/2)^m=o(1)\). For \(m>\sqrt N\), they sum to
at most \(\sum_{m>\sqrt N}2^{-m}=o(1)\). This proves \textup{(S)} with
probability tending to one.

Finally, write \(D_{\max}=\max_{i<j}D_{ij}\).
For each ordered tuple of \(\ell\) distinct vertices, the probability
that all consecutive edges and the closing edge are present is at
most \((D_{\max}/n)^\ell\).
There are at most \((rn)^\ell\) such tuples. The expected number of
cycles of length less than \(g\) is therefore at most
\(\sum_{\ell=3}^{g-1}(rD_{\max})^\ell\), a constant independent of \(n\).
Markov's inequality gives the claimed short-cycle bound with
probability tending to one.
\end{proof}

\begin{proof}[Proof of Theorem~\ref{thm:main}]
Set \(r=2d+2\). In the prescribed order, take the first parameter
\(D_*\ge\max\{D_0(r),r^3\}\) and each subsequent parameter larger
than the \(24\)th power of its predecessor. These parameters satisfy
\textup{(O)}. With the parameters fixed, choose \(n\) sufficiently
large and take a realization \(H\) as in Lemma~\ref{lem:random}.

Let \(I\) be an independent set in \(H\). If \(I\) met two blocks
\(V_i,V_j\) in at least \(n/r\) vertices each, then
\(nD_{ij}^{-1/3}\le n/r\) and \textup{(L)} would give an edge of
\(H[I]\). Thus at most one block contributes \(n/r\) or more vertices
to \(I\), and the independence number satisfies
\(\alpha(H)\le n+(r-1)n/r\le 2n\).

We use a vertex-deletion variant of the argument of
Erd\H{o}s~\cite[p.~36]{Erdos59} for constructing graphs of large girth
and chromatic number. Following Diestel~\cite[p.~301, proof of
Theorem~11.2.2]{Diestel}, choose one vertex from each cycle of length
less than \(g\) in \(H\), and let \(G\) be the induced subgraph obtained
by deleting the chosen vertices. Then \(G\) has girth at least \(g\).
At most \(\sqrt n\) vertices were deleted, and every color class of \(G\)
is independent in \(H\). Hence
\[
 \chi(G)\ge\frac{|V(G)|}{\alpha(H)}
 \ge\frac{rn-\sqrt n}{2n}\ge d.
\]
Since \(G\) is induced, any anticomplete pair in \(G\) induces the
same two subgraphs in \(H\) and is anticomplete there as well.
Lemma~\ref{lem:order} gives the desired bound.
\end{proof}

\section*{Use of generative AI}
We used OpenAI's GPT-6.2 Astra to help construct the counterexample
and develop its proof. The authors take responsibility for the
contents of the manuscript.

\end{document}